\documentclass[a4paper,11pt]{article}
\usepackage{latexsym,amssymb,enumerate,amsmath,epsfig,amsthm,dsfont}
\usepackage[margin=1in]{geometry}
\usepackage{setspace,color}
\usepackage{tikz}
\usepackage{floatrow}
\usepackage{graphicx,subfigure}
\usepackage[ruled]{algorithm2e}
\usepackage{epstopdf}
\usepackage[multidot]{grffile}
\usepackage{comment}

\newcommand{\bu}{\mathbf{u}}
\newcommand{\bp}{\mathbf{p}}
\newcommand{\bq}{\mathbf{q}}

\newtheorem{thm}{Theorem}[section]

\begin{document}
\title{An Alternating Direction Explicit Method for Time Evolution Equations with Applications to Fractional Differential Equations}
\author{Hao Liu\thanks{School of Mathematics, Georgia Institute of Technology, 686 Cherry Street, Atlanta, GA 30332-0160, USA (Email:{\bf hao.liu@math.gatech.edu}).}
\and
Shingyu Leung\thanks{Department of Mathematics, the Hong Kong University of Science and Technology, Clear Water Bay, Hong Kong (Email:{\bf masyleung@ust.hk}).}
}
\date{
	\textit{Dedicated to Professor Roland Glowinski on the occasion of his 80th birthday}
}

\maketitle

\begin{abstract}
We derive and analyze the alternating direction explicit (ADE) method for time evolution equations with the time-dependent Dirichlet boundary condition and with the zero Neumann boundary condition. The original ADE method is an additive operator splitting (AOS) method, which has been developed for treating a wide range of linear and nonlinear time evolution equations with the zero Dirichlet boundary condition. For linear equations, it has been shown to achieve the second order accuracy in time yet is unconditionally stable for an arbitrary time step size. For the boundary conditions considered in this work, we carefully construct the updating formula at grid points near the boundary of the computational domain and show that these formulas maintain the desired accuracy and the property of unconditional stability. We also construct numerical methods based on the ADE scheme for two classes of fractional differential equations. We will give numerical examples to demonstrate the simplicity and the computational efficiency of the method.
\end{abstract}

\section{Introduction}
Time evolution equations are widely used in modeling physical phenomenon and chemical reactions. It is therefore important to develop efficient and accurate numerical schemes to approximate the solution. This paper considers the approximation of the evolution of some quantity $u(\mathbf{x},t)$ on a bounded domain $\Omega$ satisfying the following initial value problem
\begin{equation}
\begin{cases}
\frac{\partial u}{\partial t}=\mathcal{L}(u,t) \, ,\\
u(\mathbf{x},0)=u_0(\mathbf{x}).
\end{cases}
\label{eq.ivp}
\end{equation}
There are various finite difference schemes developed according to the choice of the differential operator $\mathcal{L}$ \cite{str89}. The simplest approach is the explicit method given by
$$
u^{n+1}_i=u^n_i+\Delta t \mathcal{L}(u^n,t^n)
$$
with $u^n_i=u(\mathbf{x}_i,t^n)$. This scheme is developed based on the classical forward Euler time marching and therefore is at best first order accurate in time. The accuracy in space depends on the space discretization scheme used to approximate the differential operator $\mathcal{L}(u,t)$. Relating to the stability of the method, these methods usually come with a strict restriction on the time step chosen based on the mesh resolution. If the differential operator $\mathcal{L}$ is complicated, the time step used in the numerical scheme needs to be extremely small to maintain the numerical stability.

To relax the stability condition or to obtain a higher order accurate numerical solution, it is usually more favorable to consider the implicit methods or semi-implicit methods. If $\mathcal{L}$ is linear, we can easily find the inverse operator, $(I-\Delta t\mathcal{L})^{-1}$, and can directly update $u^{n+1}$ implicitly by $u^{n+1}=(I-\Delta t\mathcal{L})^{-1}(u^n)$. This implicit scheme is first order accurate in time but unconditionally stable. For a nonlinear operator $\mathcal{L}$, however, it is usually not that straightforward how to obtain an unconditionally stable numerical scheme. One possibility is to follow the regularization approach as discussed in \cite{smereka2003semi,youleu15}. A Laplacian term has been introduced in \cite{smereka2003semi} as an extra regularization for the initial value problem
$$
\begin{cases}
\frac{\partial u}{\partial t}=\alpha\Delta u-\alpha\Delta u + \mathcal{L}(u,t)\, ,\\
u(\mathbf{x},0)=v(\mathbf{x}).
\end{cases}
$$
The numerical scheme is
$$
\frac{ u^{n+1}-u^n}{\Delta t}=\alpha\Delta u^{n+1}-\alpha\Delta u^n+\mathcal{L}(u^n,t^n) \, .
$$
This scheme is first order in time and is unconditionally stable. But the choice of $\alpha$ is in general an open problem. Following a similar idea, \cite{youleu15} has introduced an extra curvature term in the level set equation for fast evolution of some curvature dependent flows. We refer all readers to the references and thereafter for a more detailed discussion of the approach.

Other than the regularization technique, another possible approach is to linearize the nonlinear operator $\mathcal{L}$ by $L(u^n)u$ and obtains the solution at the time level $t=t_{n+1}$. The Crank-Nicolson scheme is a famous approach to this linearized equation and the method is written as
$$
u^{n+1}=\left( I-\frac{\Delta t}{2}L\right)^{-1}\left( I+\frac{\Delta t}{2}L\right)u^n \, .
$$
This scheme, as the average of the implicit scheme and the explicit scheme, is second order accurate in time and unconditionally stable (at least to linear equations). But in high dimensions, it might not be trivial to invert such a matrix $\left( I-\frac{\Delta t}{2}L\right)$. Roughly speaking, there are two general operator splitting approaches developed for inverting the operator for a wide range of applications \cite{glolet89,gloleuqia15,gloleuqia16,glooshyin16,gloleuqia17}. The first group is the multiplicative operator splitting (MOS) methods which splits the complicated operator into a product of some simple operators so that each of these subproblems can be easily solved. Mathematically, we have
$$
u^{n+1}=\prod_{k=1}^K \left[ I-\Delta t L_k(u^n) \right]^{-1} u^n \, .
$$
One famous example of the MOS methods is the alternating direction implicit scheme (ADI). For example, when $Lu=\Delta u$, the ADI scheme reads as
\begin{equation*}
\begin{cases}
\frac{2}{\Delta t} \left( u^{n+\frac{1}{2}}-u^n \right) =\delta_{xx}u^{n+\frac{1}{2}}+\delta_{yy}u^n \, ,\\
\frac{2}{\Delta t} \left( u^{n+1}-u^{n+\frac{1}{2}} \right)=\delta_{xx}u^{n}+\delta_{yy}u^{n+\frac{1}{2}}
\end{cases}
\end{equation*}
where $\delta_{xx}$ and $\delta_{yy}$ denotes the second derivative operators in the $x$- and the $y$-directions, respectively. The ADI method is efficient since the method inverts only a tridiagonal matrix in each sub-step. For two dimensional problems, the ADI scheme is unconditionally stable. But for higher dimensional problems, the trivial implementation becomes only conditionally stable. Special strategy is needed to keep the unconditional stability.

Another group of splitting method is the additive operator splitting (AOS) methods \cite{luneitai91,luneitai92,vab13}. Instead of approximating the inverse by a product of simpler operators, the AOS methods split the operator into a summation of some simple operators. Mathematically, we can approximate the time evolution as
$$
u^{n+1}=\frac{1}{m} \sum_{k=1}^K \left[ I-\Delta t L_k(u^n) \right]^{-1} u^n \, .
$$
In general, the choice for the splitting might be equation-specific or might not be easy to generalize in high dimensions. One general approach is the alternating direction explicit (ADE) method. The ADE scheme was first introduced and analyzed in \cite{larkin1964some,barakat1966solution} for the linear heat equation. It has recently been extended to a much wider class of nonlinear PDEs in \cite{leuosh04}. For linear problems, the ADE method can be proved to be second order accurate in both time and space and is unconditionally stable for arbitrary time step size for both low dimensional and high dimensional problems. Using the sweeping strategy in updating the solution sequentially, the ADE scheme can be implemented in a fully explicitly fashion. In \cite{leuosh04}, the method has been carefully analyzed in the operator form and has been extended to various nonlinear equations including the fourth order diffusion equation, the Hamilton-Jacobi equation, the fourth order nonlinear lubrication equation, and etc. More recently, the method has been coupled with the fast Huygens sweeping method \cite{leuqiaser13} for nematic liquid crystal modeling \cite{klwq17}. These discussions, however, concentrate only on the constant Dirichlet boundary condition. In this work, we analyze the ADE method for the heat equation with the time dependent Dirichlet boundary condition and also the Neumann boundary condition. To demonstrate the effectiveness of the approach, we then apply the ADE scheme to solve two kinds of fractional differential equations.

This paper is organized as follows. In Section \ref{sec.ade}, we first show some existing results on the ADE method for problems with the constant Dirichlet boundary condition. Then we derive and perform a truncation error analysis of the ADE method for the heat equation with the time-dependent Dirichlet boundary condition and the Neumann boundary condition. In Section \ref{sec.tdw}, we apply our ADE method to two kinds of fractional differential equations: the time distributed order super-diffusive equation and the sub-diffusive diffusion-reaction system. Results of numerical experiments are shown in Section \ref{sec.numerical}. Finally in Section \ref{sec.conclusion} we summarize and conclude the ADE method and our results.

\section{The alternating direction explicit (ADE) method}
\label{sec.ade}

In this section, we first give a brief introduction to the ADE scheme for problem with the zero Dirichlet boundary condition. For a more detailed discussion on the numerical scheme, we refer all readers to \cite{leuosh04}. Then, we will extend the formulation to problems with the time-dependent Dirichlet boundary condition and also the zero Neumann boundary condition.

\subsection{The ADE scheme for linear time dependent PDEs}
We consider the following linear time dependent PDE
$$
\frac{\partial u}{\partial t}=\mathcal{L}(u) \, ,
$$
one first discretizes the right hand side to get
$$
\frac{\partial \bu}{\partial t}=A\bu+b
$$
where $A$ is the discretization of the linear differential operator $\mathcal{L}$, the vector $\bu$ represents the discretization of the function $u$ in space and $b$ is a constant vector. We then decompose $A$ into the form
$$
A=L+D+U
$$
where $L,D$ and $U$ are the strictly lower part, diagonal part and strictly upper part of $A$ respectively. Defining
$$
B=L+\frac{1}{2}D \, \mbox{ and } \, C=U+\frac{1}{2}D \, ,
$$
we can write the ADE scheme as the average of
\begin{equation}
(I-\Delta t \, B) \, \bp^{n+1}= \bu^n+\Delta t \,C \bu^n+\Delta t\, b
\label{eq.adesub.1}
\end{equation}
and
\begin{equation}
(I-\Delta t \, C) \,\bq^{n+1}=\bu^n+\Delta t \,B \bu^{n}+\Delta t \,b \, .
\label{eq.adesub.2}
\end{equation}
In equations (\ref{eq.adesub.1}) and (\ref{eq.adesub.2}), the matrices $(I-\Delta tB)$ and $(I-\Delta tC)$ are simply lower triangular and upper triangular, respectively. Therefore $\bp^{n+1}$ and $\bq^{n+1}$ can be updated fully explicitly using Gauss-Seidel iterations. Since the ADE method updates $u^{n+1}$ by taking the average of solutions by (\ref{eq.adesub.1}) and (\ref{eq.adesub.2}), the updating formula of the ADE scheme in operator form can be written as
\begin{eqnarray}
\bu^{n+1}&=&\frac{1}{2}\left[ \left( I-\Delta tB\right)^{-1}\left(I+\Delta tC\right) + \left( I-\Delta tC\right)^{-1} \left(I+\Delta tB\right)\right] \bu^n \nonumber\\
&&\hspace{1cm}+ \frac{1}{2}\left[ \left( I-\Delta tB\right)^{-1} + \left( I-\Delta tC\right)^{-1} \right] b \, .
\label{eq.ade.operator}
\end{eqnarray}
Because of the triangular structure of the matrices $B$ and $C$, the inverse of these operators can be easily determined using the forward substitution and backward substitution. The overall algorithm is shown in Algorithm \ref{Alg:ADEHeatZeroDirichlet} for completeness.

\begin{algorithm}
 Initialization: $n=0$ and $u^0$\;
 \While{$n<N$}{
$n \leftarrow n+1$ \;
Set $p_i^n=u_i^n$ and $q_i^n=u_i^n$ for $i=0,1,\cdots,M$ where $u_0^n=u_M^n=0$ \;
Set $p_0^{n+1}=p_M^{n+1}=q_0^{n+1}=q_M^{n+1}=0$ \;
For $i=1,2,\cdots,M-1$, solve
$$
\frac{p_i^{n+1}-p_i^n}{\Delta t}=\frac{p_{i-1}^{n+1}-p_i^{n+1}-p_i^n+p_{i+1}^n}{\Delta x^2}+\tilde{b}_i^{n+\frac{1}{2}} \, . \;
$$
For $i=M-1,M-2,\cdots,1$, solve
$$
\frac{q_i^{n+1}-q_i^n}{\Delta t}=\frac{q_{i-1}^{n}-q_i^{n}-q_i^{n+1}+q_{i+1}^{n+1}}{\Delta x^2}+\tilde{b}_i^{n+\frac{1}{2}} \, . \;
$$
For $i=M-1,M-2,\cdots,1$, compute $u_i^{n+1}=\frac{1}{2} \left( p_i^{n+1}+q_i^{n+1} \right)$.
}
 \caption{The ADE scheme for the general heat equation with the zero Dirichlet boundary condition.}
 \label{Alg:ADEHeatZeroDirichlet}
\end{algorithm}

Even though such algorithm was first introduced only for the linear heat equation, the method is applicable to a much wider range of equations such as the first order nonlinear Hamilton-Jacobi equation, a fourth order nonlinear equation and curvature dependent flows \cite{leuosh04}. The numerical scheme has various properties concerning the numerical stability and accuracy. We simply quote the following theorems from \cite{leuosh04}. For the detailed technical analysis, we refer readers to the reference thereafter.

\begin{thm} \cite{leuosh04} The ADE scheme has the following properties.
\begin{enumerate}
\item If the diagonal elements of $A$ are all non-positive, the ADE scheme (\ref{eq.ade.operator}) is second order accurate in time.
\item If $A$ is symmetric negative definite, the ADE scheme (\ref{eq.ade.operator}) is unconditionally stable.
\item If $A$ is lower-triangular with all diagonal elements negative, the ADE scheme (\ref{eq.ade.operator}) is unconditionally stable.
\end{enumerate}
\end{thm}

\subsection{The ADE scheme for the heat equation with a time dependent Dirichlet boundary condition}

In the original derivation of the numerical scheme, one considers only the time-independent Dirichlet boundary condition such that the solution is fixed (to be zero) on the boundary at all time. In this subsection, we derive the ADE method for a general heat equation in one dimension given by
\begin{equation}
\begin{cases}
\frac{\partial u}{\partial t}=\frac{\partial^2 u}{\partial x^2}+b(t) \, ,\\
u(0,t)=f(t), u(l,t)=g(t) \, ,\\
u(x,0)=v(x) \, .
\end{cases}
\label{eq.heat}
\end{equation}
In equation (\ref{eq.heat}), both the source term and the boundary conditions are now functions of $t$. The corresponding derivation to higher dimensions is rather straightforward and will not be discussed in the current work. Let $\Omega$ be discretized using $x_i$'s for $i=0,1,\cdots,M$ such that $x_i=i\Delta x$, $x_0=0$ and $x_M=l$. Using a time step $\Delta t$, we denote $t_n=n\Delta t$, for $n=0,1,\cdots,N$ such that $t_N=T$. Furthermore, we introduce the following notations given by
\begin{eqnarray*}
u_{i}^n=u(x_i,t_n) \, , \, \ u_{i}^{n+\frac{1}{2}}=\frac{1}{2}(u_{i}^n+u_{i}^{n+1}) \, \mbox{ and } \, \tilde{u}_{i}^{n+\frac{1}{2}}=u\left(x_i,\frac{1}{2}(t_{n}+t^{n+1})\right) \, .
\end{eqnarray*}
Then, equation (\ref{eq.heat}) can be discretized using the ADE scheme approach given by
$$
\begin{cases}
\frac{u^{n+1}-u^n}{\Delta t}=ADE(\Delta_h,u^n,u^{n+1})+\tilde{b}^{n+\frac{1}{2}} \, ,\\
u_0^n=f^n, u_M^n=g^n \, ,\\
u_i^0=v_i \, , \, i=0,1,\cdots,M \, \mbox{ and } \, n=0,1,\cdots,N \, ,
\end{cases}
$$
where the term $ADE(\Delta_h,u^n,u^{n+1})$ denotes the ADE treatment of the discrete Laplacian $\Delta_h$ involving both the directional sweepings and the averaging step which is similar to that in the linear heat equation. To summarize, the method is given in Algorithm \ref{Alg:ADEHeatTimeDirichlet}.

\begin{algorithm}
 Initialization: $n=0$ and $u^0$\;
 \While{$n<N$}{
$n \leftarrow n+1$ \;
Set $p_i^n=u_i^n$ and $q_i^n=u_i^n$ for $i=0,1,\cdots,M$ where $u_0^n=f^n$ and $u_M^n=g^n$ \;
Set $p_0^{n+1}=f^{n+1}$ and $q_M^{n+1}=g^{n+1}$ satisfying the boundary conditions at $x=x_0$ and $x_M$, respectively\;
For $i=1,2,\cdots,M-1$, solve
\begin{equation}
\frac{p_i^{n+1}-p_i^n}{\Delta t}=\frac{p_{i-1}^{n+1}-p_i^{n+1}-p_i^n+p_{i+1}^n}{\Delta x^2}+\tilde{b}_i^{n+\frac{1}{2}} \, . \;
\label{eq.ade.heat.1}
\end{equation}
For $i=M-1,M-2,\cdots,1$, solve
\begin{equation}
\frac{q_i^{n+1}-q_i^n}{\Delta t}=\frac{q_{i-1}^{n}-q_i^{n}-q_i^{n+1}+q_{i+1}^{n+1}}{\Delta x^2}+\tilde{b}_i^{n+\frac{1}{2}} \, . \;
\label{eq.ade.heat.2}
\end{equation}
For $i=M-1,M-2,\cdots,1$, compute
\begin{equation}
u_i^{n+1}=\frac{1}{2} \left( p_i^{n+1}+q_i^{n+1} \right) \, . \;
\label{eq.ade.heat.3}
\end{equation}
}
 \caption{The ADE scheme for the general heat equation with the time-dependent Dirichlet boundary condition.}
 \label{Alg:ADEHeatTimeDirichlet}
\end{algorithm}

The scheme (\ref{eq.ade.heat.1})-(\ref{eq.ade.heat.3}) is similar to the original ADE method for the heat equation with the time-independent Dirichlet boundary condition. The main difference is only in the treatment to those grid points adjacent to the boundary. In equation (\ref{eq.ade.heat.1}), $p_i^n$ is updated in the ascending direction from $i=1$ to $i=M-1$. For $i=1$, since we have the boundary condition $u_0^{n+1}=f^{n+1}$, equation (\ref{eq.ade.heat.1}) can be written as
\begin{equation}
\frac{p_1^{n+1}-p_1^n}{\Delta t}=\frac{f^{n+1}-p_1^{n+1}-p_1^n+p_2^n}{\Delta x^2}+\tilde{b}_i^{n+\frac{1}{2}} \, .
\label{eq.ade.heat.1.1}
\end{equation}
Therefore, the quantity $p_1^{n+1}$ is updated using the exact boundary condition $f^{n+1}$ at $t=t_{n+1}$. For the grid point with the index $i=M-1$, we use the boundary condition $u_M^n=g^n$ and write equation (\ref{eq.ade.heat.1}) as
\begin{equation}
\frac{p_{M-1}^{n+1}-p_{M-1}^n}{\Delta t}=\frac{p_{M-2}^{n+1}-p_{M-1}^{n+1}-p_{M-1}^n+g^n}{\Delta x^2}+\tilde{b}_i^{n+\frac{1}{2}} \, .
\label{eq.ade.heat.1.2}
\end{equation}
Therefore, the quantity $p_{M-1}^n$ is updated using the value of $g$ at the previous time level given by $g^n$ at $t=t_n$. Although we have the exact boundary condition at $t=t_{n+1}$, we use its value at $t=t_n$.

Similar consideration applies to the update of $q$ according to equation (\ref{eq.ade.heat.2}) where we sweep through the index in a descending direction from $i=M-1$ to $i=1$. For the index $i=M-1$, we have the boundary condition $u_M^{n+1}=g^{n+1}$ so that equation (\ref{eq.ade.heat.2}) can be written as
\begin{equation}
\frac{q_{M-1}^{n+1}-q_{M-1}^n}{\Delta t}=\frac{q_{M-2}^{n}-q_{M-1}^{n}-q_{M-1}^{n+1}+g^{n+1}}{\Delta x^2}+\tilde{b}_i^{n+\frac{1}{2}} \, .
\label{eq.ade.heat.2.1}
\end{equation}
The quantity $q_{M-1}^{n+1}$, therefore, is updated using the exact boundary condition at the new time level $t=t_{n+1}$. For the left end point at the index $i=1$, we impose the boundary condition $u_0^n=f^n$ for equation (\ref{eq.ade.heat.2}) and rewrite it as
\begin{equation}
\frac{q_1^{n+1}-q_1^n}{\Delta t}=\frac{f^{n}-q_1^{n}-q_1^{n+1}+q_{2}^{n+1}}{\Delta x^2}+\tilde{b}_i^{n+\frac{1}{2}} \, .
\label{eq.ade.heat.2.2}
\end{equation}
This implies that the quantity $q_1^{n+1}$ is updated using the boundary condition imposed at the previous time level $t=t_n$.

Now we analyze the numerical accuracy of the proposed ADE scheme with the time-dependent boundary condition. One main property of this scheme is that the above updating formulas (\ref{eq.ade.heat.1.1})-(\ref{eq.ade.heat.2.2}) at grid points adjacent to the boundary are designed to keep the symmetry of scheme (\ref{eq.ade.heat.1})-(\ref{eq.ade.heat.3}). By Taylor expansion at $(x=i\Delta x, t=(n+\frac{1}{2})\Delta t)$, the leading order truncation error in $p$ and $q$ are given by
$$
-\frac{\Delta t}{\Delta x}p_{tx} -\frac{1}{6}\Delta t \Delta x p_{txxx}-\frac{1}{24}\Delta t^2\left( \frac{\Delta t}{\Delta x}\right) p_{tttx}+\Delta t^2\left(\frac{1}{8} p_{ttxx}-\frac{1}{24}p_{ttt}\right)+\frac{1}{12}\Delta x^2 p_{xxxx}
$$
and
$$
\frac{\Delta t}{\Delta x}q_{tx} +\frac{1}{6}\Delta t \Delta x q_{txxx}+\frac{1}{24}\Delta t^2\left( \frac{\Delta t}{\Delta x}\right) q_{tttx}+\Delta t^2\left(\frac{1}{8} q_{ttxx}-\frac{1}{24}q_{ttt}\right)+\frac{1}{12}\Delta x^2 q_{xxxx} \, ,
$$
respectively. After the averaging step (\ref{eq.ade.heat.3}), the leading error of the overall scheme (\ref{eq.ade.heat.1})-(\ref{eq.ade.heat.3}) is therefore given by
$$
\Delta t^2\left(\frac{1}{8} u_{ttxx}-\frac{1}{24}u_{ttt}\right)+\frac{1}{12}\Delta x^2 u_{xxxx} \, .
$$


In the above derivation, we assume that the function $b$ can be evaluated at all half-time level given by $\tilde{b}_{i}^{n+\frac{1}{2}}$. In some applications, however, when we can only provide $b$ at gridded levels $t=n\Delta t$ for some integer $n$, we approximate the quantity by $b_i^{n+\frac{1}{2}}$ which is second order in time. Then the leading error of the scheme (\ref{eq.ade.heat.1})-(\ref{eq.ade.heat.3}) is therefore given by
$$
\Delta t^2\left(\frac{1}{16} u_{ttxx}-\frac{1}{24}u_{ttt}+\frac{1}{4}b_{tt}\right)+\frac{1}{12}\Delta x^2 u_{xxxx} \, ,
$$
which is still $O(\Delta t^2,\Delta x^2)$ accurate.

Now, we discuss the overall computational complexity. Since the updating formulas involve only explicit operation at each individual grid point, the computational complexity of obtaining $u_i^{n+1}$ for all $i$ from $u_i^{n}$ is of $O(M)$ with the total number of grid points is given by $M+1$. This means that the complexity is in fact the same as typical explicit scheme for the heat equation. But because the method is unconditionally stable, it takes less time steps to reach the final time of the simulation.

\subsection{The ADE scheme for the heat equation with the zero Neumann boundary condition}
In this subsection, we discuss the implementation of the Neumann boundary condition. Even though we consider only the zero Neumann boundary condition, the algorithm can be easily modified to solve problems with non-zero Neumann boundary condition.

We consider the same problem given by equation (\ref{eq.heat}) but replace the Dirichlet boundary condition by the zero Neumann conditions given by $u_x(0)=0$ and $u_x(l)=0$. We further assume that we have already obtained $u^n$ and would like to use the ADE scheme to get $u^{n+1}$. Numerically, the zero Neumann boundary condition is typically approximated by standard finite difference and is given by $u_1^n=u_0^n$ and $u_M^n=u_{M-1}^n$ for all $n$. In the first sub-problem of the ADE scheme where we sweep through the solution in the ascending order, we need the boundary condition $p_0^{n+1}$ at the updated time level $t_{n+1}$. We note that such value is actually an unknown value if we require the solution satisfies $p_0^{n+1}=p_1^{n+1}$. However, if we substitute this constraint to the updating formula of the ADE scheme for the first sub-problem, we can actually obtain
\begin{equation}
p_0^{n+1} = \left( 1 -\frac{\Delta t}{\Delta x^2} \right) p_1^n+ \frac{\Delta t}{\Delta x^2} \, p_2^n+ \Delta t \, \tilde{b}_1^{n+\frac{1}{2}} \, .
\label{Eqn:NeumannBCLeft}
\end{equation}
Such choice of the boundary condition will automatically lead to the zero Neumann boundary condition at the left boundary. Similarly, for the second sub-problem when we sweep through the index in the descending order, we impose the following boundary condition at the right boundary given by
\begin{equation}
q_M^{n+1} = \frac{\Delta t}{\Delta x^2} \, q_{M-2}^n+ \left(1-\frac{\Delta t}{\Delta x^2} \right) q_{M-1}^n+ \Delta t \, \tilde{b}_M^{n+\frac{1}{2}} \, .
\label{Eqn:NeumannBCRight}
\end{equation}
Such choice will automatically lead to $q_M^{n+1}=q_{M-1}^{n+1}$. To summarize, we give the overall scheme in Algorithm \ref{Alg:ADEHeatNeumann}. The scheme is extremely easy to implement. Comparing to Algorithm \ref{Alg:ADEHeatTimeDirichlet}, we only need to modify boundary conditions by imposing $u_0^n=u_1^n$ and $u_M^n=u_{M-1}^n$, setting $p_0^{n+1}$ using (\ref{Eqn:NeumannBCLeft}) and $q_M^{n+1}$ using (\ref{Eqn:NeumannBCRight}). The main part of the numerical implementation remains exactly the same.

\begin{algorithm}
 Initialization: $n=0$ and $u^0$\;
 \While{$n<N$}{
$n \leftarrow n+1$ \;
Set $p_i^n=u_i^n$ and $q_i^n=u_i^n$ for $i=0,1,\cdots,M$ where $u_0^n=u_1^n$ and $u_M^n=u_{M-1}^n$ \;
Set $p_0^{n+1}$ using (\ref{Eqn:NeumannBCLeft}) and $q_M^{n+1}$ using (\ref{Eqn:NeumannBCRight}) \;
For $i=1,2,\cdots,M-1$, solve
$$
\frac{p_i^{n+1}-p_i^n}{\Delta t}=\frac{p_{i-1}^{n+1}-p_i^{n+1}-p_i^n+p_{i+1}^n}{\Delta x^2}+\tilde{b}_i^{n+\frac{1}{2}} \, . \;
$$
For $i=M-1,M-2,\cdots,1$, solve
$$
\frac{q_i^{n+1}-q_i^n}{\Delta t}=\frac{q_{i-1}^{n}-q_i^{n}-q_i^{n+1}+q_{i+1}^{n+1}}{\Delta x^2}+\tilde{b}_i^{n+\frac{1}{2}} \, . \;
$$
For $i=M-1,M-2,\cdots,1$, compute $u_i^{n+1}=\frac{1}{2} \left(p_i^{n+1}+q_i^{n+1} \right)$.
}
 \caption{The ADE scheme for the general heat equation with the zero Neumann boundary condition.}
 \label{Alg:ADEHeatNeumann}
\end{algorithm}

\section{Applications to partial differential equations with fractional derivatives}

Because of the nice stability property of the ADE scheme, the method is especially suitable to problems for which standard finite difference methods require a restrictive time stability constraint. In this section, we discuss two applications of the ADE scheme related to fractional derivatives. The first one is a time-distributed order super-diffusive partial differential equation with the time dependent Dirichlet boundary condition. The second one is a reaction-diffusion system in a sub-diffusive regime with the zero Neumann boundary condition.

\subsection{The time distributed order super-diffusive partial differential equation}
\label{sec.tdw}

Fractional differential equation (FDE) is a powerful technique in modeling an anomalous phenomena such as anomalous diffusion and anomalous transportation \cite{milros93,wys96,pod98}. In the diffusion process, the mean squared displacement (MSD) can be written as an exponential function of the time $t$ with a coefficient $\alpha$. For the case where $\alpha=1$, the diffusion process reduces back to the typical isotropic diffusion modeled by the heat equation $u_t=\Delta u$. The anomalous diffusion includes two types: sub-diffusion and super-diffusion. The diffusion process is called sub-diffusion (diffuse slower than normal) if $0<\alpha<1$, and is called super-diffusion (diffuse faster than normal) if $1<\alpha<2$. Mathematically, these anomalous diffusion processes are modeled using a time fractional order differential equation \cite{eab2011fractional}. More recently, some efforts are spent to solve more complicated problems where the diffusion process contains diffusions in various time scales so that the diffusion speed is modeled by a range of parameter $\alpha$. This leads to the so-called time distributed order differential equations \cite{gao2016two}. In this section, we consider the following time distributed-order super-diffusive problem given by
\begin{equation}
\begin{cases}
\int_1^2 w(\gamma) {}_0^C\!D_t^{\gamma} u(x,y,t)d\gamma=\Delta u(x,y,t)+F(x,y,t), (x,y)\in\Omega,0<t\leq T,\\
u(x,y,0)=0, u_t(x,y,0)=0, (x,y)\in \Omega \, ,\\
u(x,y,t)=\psi(x,y,t), (x,y)\mbox{ on }\partial\Omega, 0\leq t\leq T \, ,
\end{cases}
\label{eq.fde}
\end{equation}
with $\Omega=\{(x,y)| 0<x<L_1,0<y<L_2\}$, and $\partial\Omega$ being the boundary of $\Omega$. The functions $F(x,y,t)$ and $\psi(x,y,t)$ are given to model the time-dependent source term and the initial condition, respectively. The function $w(\gamma)\geq0$ with $\int_1^2 w(\gamma)d\gamma=c_0>0$ determines the weights among various diffusion speed scales. In equation (\ref{eq.fde}), the derivative ${}_0^C\!D_t^{\gamma} u(x,y,t)$ is the $\gamma$-th order time Caputo fractional derivative defined as
\begin{equation}
{}_0^C\!D_t^{\gamma} u(x,y,t)=
\begin{cases}
u_t(x,y,t)-u_t(x,y,0), \mbox{ if } \gamma=1 \, ,\\
\frac{1}{\Gamma(2-\gamma)} \int_0^t (t-\xi)^{1-\gamma} \frac{\partial^2 u}{\partial \xi^2}(x,y,\xi) d\xi, 1<\gamma<2,\\
u_{tt}(x,y,t), \mbox{ if } \gamma=2 \, .
\end{cases}
\end{equation}

Consider a rectangle domain $\Omega=\{(x,y) | 0\leq x\leq L_1,0\leq y\leq L_2\}$. Let $\Delta x$ and $\Delta y$ be the mesh size in the $x$- and the $y$-direction, respectively. Define $x_i=i\Delta x$ and $y_j=j\Delta y$ for $0\leq i\leq M_1$ and $0\leq j\leq M_2$ such that $x_{M_1}=L_1$ and $y_{M_2}=L_2$. Denote the set $\omega=\{ (i,j) |1\leq i \leq M_1-1, 1\leq j \leq M_2-1\}$ and $\omega_0=\{(i,j)| (x_i ,y_j)\in \partial\Omega\}$ and $\omega_h=\omega\cap\omega_0$. For any $v\in \Omega_h$, the discretization of $\Omega$, we introduce the following difference notations:
\begin{eqnarray*}
&\delta_x v_{i-\frac{1}{2}j}=\frac{1}{\Delta x}(v_{ij}-v_{i-1j}),\  & \delta_x^2 v_{ij}=\frac{1}{\Delta x}\left( \delta_x v_{i+\frac{1}{2}j}- \delta_x v_{i-\frac{1}{2}j} \right),\\
&\delta_y v_{i-\frac{1}{2}j}=\frac{1}{\Delta y}(v_{ij}-v_{ij-1}),\  & \delta_y^2 v_{ij}=\frac{1}{\Delta y}\left( \delta_x v_{ij+\frac{1}{2}}- \delta_y v_{ij-\frac{1}{2}} \right),\\
&\Delta_h v_{ij}=\delta_x^2 v_{ij}+ \delta_y^2 v_{ij}\, .&
\end{eqnarray*}

In the time space, we use $\Delta t$ to denote our time step. Define $t_n=n\Delta t, n=0,1,\cdots,N$ such that $t_N=T$. For the integral on the fractional derivative in time, we divide the integration domain $[1,2]$ into $J$ subintervals with length $\Delta \gamma=\frac{1}{J}$ and define $\gamma_l=1+l\Delta \gamma, l=0,1,\cdots,J.$ Then the integration will be approximated using Trapezoidal rule:
\begin{equation*}
\int_1^2 s(\gamma)d\gamma=\Delta\gamma \sum_{l=0}^{J} c_l s_{\gamma_l}+O(\Delta \gamma^2)
\end{equation*}
with $c_l=\frac{1}{2}$ for $l=0$ or $J$, and $c_l=1$ otherwise.

Following a similar strategy as in \cite{gao2016two}, for $\alpha\in\left[0,1\right]$, we define $g_k^{(\alpha)}$, $\lambda_k^{(\alpha)}$ and $\mu$ as
\begin{eqnarray*}
&&g_0^{(0)}=1, \ g_k^{(0)}=0  \, , \,  g_0^{\alpha}=1,\ g_k^{\alpha}=\left(1-\frac{\alpha+1}{k}\right) g_{k-1}^{(\alpha)} \, , \, \lambda_0^{(\lambda)}=\left( 1+\frac{\alpha}{2}\right) g_0^{(\alpha)} \, ,\\
& &\lambda_k^{(\alpha)}= \left( 1+\frac{\alpha}{2}\right) g_k^{(\alpha)} -\frac{\alpha}{2} g_{k-1}^{(\alpha)} \, , \, \mu= \Delta \gamma \sum_{l=0}^{J} c_l w(\gamma_l) \tau^{-\alpha_l} \lambda_0^{(\alpha_l)} \, ,
\end{eqnarray*}
for $k\geq1$ and $0<\alpha\leq1$. Denote $\alpha_l=\gamma_l-1$, a simple finite difference scheme applied to equation (\ref{eq.fde}) is given by
$$
\begin{cases}
\Delta\gamma \sum_{l=0}^{J} c_l w(\gamma_l)  \Delta t^{-\alpha_l}  \sum_{k=0}^{n} \lambda_k^{(\alpha_l)} \delta_t u_{ij}^{n-k+\frac{1}{2}} = \Delta_h u_{ij}^{n+\frac{1}{2}} +F_{ij}^{n+\frac{1}{2}}, (i,j)\in \omega, 1\leq n \leq N,\\
u_{ij}^0=0, (i,j)\in \omega \, \mbox{ and } \, u_{ij}^n=\psi(x_i,y_j,t_n), (i,j)\in \omega_0, 0\leq n \leq N-1 \, ,
\end{cases}
$$
where $\Delta_h u_{ij}^{n+\frac{1}{2}}$ is an approximation of the Laplacian of the solution at the half time level $t_{n+\frac{1}{2}}$ using methods such as central difference. The main issue, however, is that it is not easy to find an efficient solver to invert for $u^{n+1}$. In particular, ADI scheme is not directly applicable. In a recent work, \cite{gao2016two} has introduced an extra regularization to the discretization by $\frac{\Delta t^2}{4\mu} \delta_x^2\delta_y^2 \delta_t u_{ij}^{n+\frac{1}{2}}$ which also nicely coupled with the ADI scheme,
\begin{equation}
\begin{cases}
{\displaystyle \Delta\gamma \sum_{l=0}^{J} c_l w(\gamma_l) \Delta t^{-\alpha_l} \sum_{k=0}^{n} \lambda_k^{(\alpha_l)} \delta_t u_{ij}^{n-k+\frac{1}{2}} +\frac{\Delta t^2}{4\mu} \delta_x^2\delta_y^2 \delta_t u_{ij}^{n+\frac{1}{2}} = \Delta_h u_{ij}^{n+\frac{1}{2}} +F_{ij}^{n+\frac{1}{2}} \, ,}\\
u_{ij}^0=0, (i,j)\in \omega \, \mbox{ and } \, u_{ij}^n=\psi(x_i,y_j,t_n), (i,j)\in \omega_0, 0\leq n \leq N-1 \,.
\end{cases}
\label{eq.adi}
\end{equation}
This ADI scheme (\ref{eq.adi}) as developed in \cite{gao2016two} has been proven to give an approximation to (\ref{eq.fde}) with the error $O(\Delta t^2+\Delta x^2 +\Delta y^2 +\Delta \gamma^2)$. In the derivation of (\ref{eq.adi}), the extra regularization term $\frac{\Delta t^2}{4\mu} \delta_x^2\delta_y^2 \delta_t u_{ij}^{n+\frac{1}{2}}$ is manually incorporated into the equation. Such term is necessary in the ADI scheme solely originated from the construction of the numerical method but is unnecessary in our ADE approach.

Our proposed ADE finite difference scheme is given by
\begin{equation}
\begin{cases}
{\displaystyle \Delta\gamma \sum_{l=0}^{J} c_l w(\gamma_l) \Delta t^{-\alpha_l} \sum_{k=0}^{n} \lambda_k^{(\alpha_l)} \delta_t u_{ij}^{n-k+\frac{1}{2}} = ADE(\Delta_h,u^{n},u^{n+1}) +F_{ij}^{n+\frac{1}{2}} \, ,}\\
u_{ij}^0=0, (i,j)\in \omega \, \mbox{ and } \, u_{ij}^n=\psi(x_i,y_j,t_n), (i,j)\in \omega_0, 0\leq n \leq N-1 \, .
\end{cases}
\label{eq.ade.scheme}
\end{equation}
Note that this expression can be reorganized as
\begin{equation}
\begin{cases}
{\displaystyle
\mu\delta_t u_{ij}^{n+\frac{1}{2}} = ADE(\Delta_h,u^{n},u^{n+1}) + \Delta\gamma\sum_{l=0}^{J} c_l w(\gamma_l) \Delta t^{-\alpha_l} \sum_{k=1}^{n} \lambda_k^{ (\alpha_l)} \delta_t u_{ij}^{n-k+\frac{1}{2}} +F_{ij}^{n+\frac{1}{2}} \, ,}\\
u_{ij}^0=0, (i,j)\in \omega \, \mbox{ and } \, u_{ij}^n=\psi(x_i,y_j,t_n), (i,j)\in \omega_0, 0\leq n \leq N-1
\end{cases}
\label{eq.ade.scheme1}
\end{equation}
which is a diffusion equation of $u^{n+\frac{1}{2}}$. Assume that we are given $u^k$, for $k=0,1,\cdots,n$ and want to compute $u^{n+1}$. In equation (\ref{eq.ade.scheme1}), the time derivative term $\delta_t u_{ij}^{k+\frac{1}{2}}$ can be approximated by
\begin{equation}
\delta_t u_{ij}^{k+\frac{1}{2}}=\frac{u_{ij}^{k+1}-u_{ij}^k}{\Delta t} \, ,
\label{eq.timeDis}
\end{equation}
for $k=0,1,\cdots,n,$ with the second order accuracy in time. We can analytically determine $F_{ij}^{n+\frac{1}{2}}$ since the expression of $F$ is given. Substituting the approximation (\ref{eq.timeDis}) to the equation (\ref{eq.ade.scheme1}) and applying the ADE method, we finally get the following Algorithm \ref{Alg:ADETimeDistributed}.

\begin{algorithm}
 Initialization: $n=0$ and $u^0$\;
 \While{$n<N$}{
$n \leftarrow n+1$ \;
Given $u^k$ for $k=0,1,\cdots,n$ and also the function $F$\;
Compute $\delta_t u^{n-k+\frac{1}{2}},k=1,2,\cdots,n$ and $F^{n+\frac{1}{2}}$\;
For $i=1,2,\cdots,M_1-1$ and $j=1,2,\cdots,M_2-1$, compute
$$
\tilde{b}_{ij}^{n+\frac{1}{2}}=-\Delta\gamma\sum_{l=0}^{J} c_l w(\gamma_l) \Delta t^{-\alpha_l}\sum_{k=1}^{n} \lambda_k^{ (\alpha_l)} \delta_t u_{ij}^{n-k+\frac{1}{2}} +F_{ij}^{n+\frac{1}{2}} \, .\;
$$
Set $p_{ij}^{n}=q_{ij}^{n}=v_{ij}^{n}=w_{ij}^{n}=u_{ij}^{n}$, for $i=1,2,\cdots,M_1-1$ and $j=1,2,\cdots,M_2-1$\;
Set $p_{ij}^{k}=\psi_{ij}^{k},q_{ij}^{k}=\psi_{ij}^{k}, v_{ij}^{k}=\psi_{ij}^{k},w_{ij}^{k}=\psi_{ij}^{k},\forall (i,j)\in\omega_0,k=n+1,n$\;
For $i=1,2,\cdots,M_1-1,\ j=1,2,\cdots,M_2-1$, solve
$$
\frac{p_{ij}^{n+1}-p_{ij}^{n}}{\Delta t}=\frac{p_{i-1j}^{n+1}-p_{ij}^{n+1}-p_{ij}^{n}+p_{i+1j}^{n}}{\Delta x^2} +\frac{p_{ij-1}^{n+1}-p_{ij}^{n+1}-p_{ij}^{n}+p_{ij+1}^{n}}{\Delta y^2} +\tilde{b}_i^{n+\frac{1}{2}} \, .\;
$$
For $i=1,2,\cdots,M_1-1,\ j=M_2-1,M_2-2,\cdots,1$, solve
$$
\frac{q_{ij}^{n+1}-q_{ij}^{n}}{\Delta t}=\frac{q_{i-1j}^{n+1}-q_{ij}^{n+1}-q_{ij}^{n}+q_{i+1j}^{n}}{\Delta x^2} +\frac{q_{ij-1}^{n}-q_{ij}^{n}-q_{ij}^{n+1}+q_{ij+1}^{n+1}}{\Delta y^2} +\tilde{b}_i^{n+\frac{1}{2}} \, .\;
$$
For $i=M_1-1,M_1-2,\cdots,1,\ j=1,2,\cdots,M_2-1$, solve
$$
\frac{v_{ij}^{n+1}-v_{ij}^{n}}{\Delta t}=\frac{v_{i-1j}^{n}-v_{ij}^{n}-v_{ij}^{n+1}+v_{i+1j}^{n+1}}{\Delta x^2} +\frac{v_{ij-1}^{n+1}-v_{ij}^{n+1}-v_{ij}^{n}+v_{ij+1}^{n}}{\Delta y^2} +\tilde{b}_i^{n+\frac{1}{2}} \, .\;
$$
For $i=M_1-1,M_1-2,\cdots,1,\ j=M_2-1,M_2-2,\cdots,1$, solve
$$
\frac{w_{ij}^{n+1}-w_{ij}^{n}}{\Delta t}=\frac{w_{i-1j}^{n}-w_{ij}^{n}-w_{ij}^{n+1}+w_{i+1j}^{n+1}}{\Delta x^2} +\frac{w_{ij-1}^{n}-w_{ij}^{n}-w_{ij}^{n+1}+w_{ij+1}^{n+1}}{\Delta y^2} +\tilde{b}_i^{n+\frac{1}{2}} \, .\;
$$
Compute $u_{ij}^{n+1}=\frac{1}{4} \left(p_{ij}^{n+1}+q_{ij}^{n+1}+v_{ij}^{n+1}+w_{ij}^{n+1}\right)$ for $(i,j)\in\omega$.
}
 \caption{The ADE scheme for the two dimensional time distributed order super-diffusive PDE.}
 \label{Alg:ADETimeDistributed}
\end{algorithm}

To end this section, we present two theorems of the stability and convergence of (\ref{eq.ade.scheme}). The proofs are similar to those as shown in \cite{gao2016two}. We recommend interested readers to the reference and thereafter.
\begin{thm}
  Let $\{u_{ij}^{n+1}| (i,j)\in \omega_h, 0\leq N\}$ be the solution of the following difference scheme
$$
\begin{cases}
{\displaystyle
\Delta\gamma \sum_{l=0}^{J} c_l w(\gamma_l) \Delta t^{-\alpha_l} \sum_{k=0}^{n} \lambda_k^{(\alpha_l)} \delta_t u_{ij}^{n-k+\frac{1}{2}} = ADE(\Delta_h,u^{n},u^{n+1}) +G_{ij}^{n+1} \, ,}\\
u_{ij}^0=\phi_{i,j}, (i,j)\in \omega \, \mbox{ and } \, u_{ij}^n=0, (i,j)\in \omega_0, 0\leq n \leq N-1,
\end{cases}
$$
then it holds
$$
\| \Delta_h u^{n+1}\|^2 \leq \exp(T)\left[ 3\| \Delta_h u^0\|^2 + 4(\|G^1\|^2 + \max_{1\leq k \leq n+1} \|G^k\|^2) +4\Delta \sum_{k=1}^{n} \|\delta_t G^{k+\frac{1}{2}}\|^2 \right]
$$
for $1\leq n \leq N-1$ where
$$
\| G^n\|^2=\Delta x \Delta y \sum_{i=1}^{M_1-1} \sum_{j=1}^{M_2-1} (G_{ij}^n)^2
\, \mbox{ and } \,
\| \delta_t G^{n+\frac{1}{2}} \|^2 =\Delta x \Delta y \sum_{i=1}^{M_1-1} \sum_{j=1}^{M_2-1} (\delta_t G_{ij}^{n+\frac{1}{2}})^2 \, .
$$
\end{thm}

\begin{thm}
Let $u(x,y,t)$ be the $C^2$-solution of the problem (\ref{eq.fde}) and $\{u_{ij}^n| (i,j)\in \omega_h, 0\leq n \leq N\}$ be the solution of (\ref{eq.ade.scheme}). Denote $e_{ij}^n= U_{ij}^n-u_{ij}^n, (i,j)\in \omega_h, 0\leq n \leq N$. Then we have
$$
\|e^n\|_{\infty} \leq 2C\cdot \exp(T/2) \sqrt{(2+T)L_1 L_2}(\Delta t^2+ \Delta x^2 +\Delta y^2 +\Delta \gamma^2) \, ,
$$
for $1\leq n \leq N$ where $C$ is some constant.
\end{thm}

Indeed, numerical methods for solving fractional derivative equations grow tremendously in recent years. For example, see \cite{forsim01,dffl05,lanhen05,lsat05,sunwu06,linxu07,lixu09,dddglm15}. We are not aiming to provide a complete review of the numerical approaches to the equation, nor to compare the efficiency or the accuracy of these numerical schemes. The main purpose of this section is to show that it is straightforward to implement the ADE scheme to solve this class of equations and to demonstrates the effectiveness of the scheme. We leave it as a future work to investigate the performance among various numerical schemes.

\subsection{A reaction-diffusion system in a sub-diffusive regime}

In the second application, we apply the ADE scheme to solve a sub-diffusive reaction-diffusion system for pattern formation based on diffusion-driven instability. We consider the Turing pattern \cite{turing1952chemical} which has been a very successive model to explain the patterns on various surfaces \cite{kondo1995reaction,sick2006wnt}, porous media flows \cite{drazer1999experimental}, and biological systems \cite{metzler1998multiple}. One famous class of such reaction-diffusion systems is the so-called two-component activator-inhibitor system. In this system under some conditions, one component stimulates the growth rate of both components while the other one inhibits their productions. In a recent work, \cite{hernandez2009dynamics} has generalized the standard activator-inhibitor system developed in \cite{barrio1999two} by replacing the typical diffusion by sub-diffusion. The system reads as
\begin{equation}
\begin{cases}
\frac{\partial^{\alpha}u}{\partial t^{\alpha}}=D\delta\Delta u+a_{11}u+ a_{12}v -r_2 uv-(a_{11}r_1)uv^2 \, ,\\
\frac{\partial^{\beta}v}{\partial t^{\beta}}= \delta \Delta v+ a_{21}u+ a_{22}v +r_2uv +(a_{11}r_1)uv^2 \, ,
\end{cases}
\label{eq.turing}
\end{equation}
for $0<\alpha,\beta\leq1$. The parameters $D$ and $\delta$ are two coefficients governing the diffusion associated to two components. The derivative $\frac{\partial^{\alpha}}{\partial t^{\alpha}}$ is the Caputo-type fractional derivative of order $\alpha$. For problems with finite domain, we impose the zero Neumann boundary condition on the boundary of the computational domain.

To approximate the fractional derivative numerically, we follow the idea as introduced in \cite{gorenflo2007convergence} by
\begin{equation}
\left.\frac{\partial^{\alpha}u_{ij}}{\partial t^{\alpha}} \right|_{t_{n+1}}=\sum_{m=0}^{n+1} (-1)^m\binom{\alpha}{m} \frac{u_{ij}(t_{n+1-m})-u_{ij}(t_0)}{(\Delta t)^{\alpha}} \, .
\label{eq.approxFD}
\end{equation}

In our current implementation, we treat the terms related to $v$ in the first equation in (\ref{eq.turing}) explicitly when updating the quantity $u$, while we take all $u$-related terms in the second equation in (\ref{eq.turing}) explicitly when updating $v$. Then equation (\ref{eq.turing}) can be discretized as
$$
\begin{cases}
\frac{u_{ij}^{n+1}-u_{ij}^0}{\Delta t^{\alpha}}=D\delta \, ADE(\Delta_h,u^n,u^{n+1}) +a_{11}\frac{u_{ij}^{n+1}+u_{ij}^n}{2}+ a_{12}v_{ij}^n- r_2u_{ij}^n v_{ij}^n -(a_{11}r_1)u_{ij}^n(v_{ij}^n)^2\\
\hspace{3cm}-\sum_{m=1}^n (-1)^m\binom{\alpha}{m} \frac{u_{ij}^{n+1-m}-u_{ij}^0}{(\Delta t)^{\alpha}},\\
\frac{v_{ij}^{n+1}-v_{ij}^0}{\Delta t^{\beta}}=\delta \, ADE(\Delta_h,v^n,v^{n+1}) +a_{21}u^n +a_{22}\frac{v_{ij}^{n+1}+v_{ij}^n}{2} +r_2 u_{ij}^n v_{ij}^n +(a_{11}r_1)u_{ij}^n (v_{ij}^n)^2 \\
\hspace{3cm}- \sum_{m=1}^n (-1)^m\binom{\beta}{m} \frac{v_{ij}^{n+1-m}-v_{ij}^0}{(\Delta t)^{\beta}} \, .
\end{cases}
$$
Indeed, various other ADE forms are possible for this equation in treating the nonlinear terms. But we tend to concentrate only on the effectiveness of the simple ADE discretization to complicated equation. More detailed studies on the ADE discretization of a general nonlinear system of PDEs will be given in the future.

\section{Numerical experiments}
\label{sec.numerical}
In this section, we test the accuracy of the ADE algorithm. For one dimensional examples, we use $\Delta x=1/M$ such that there are $M+1$ grid points. For two dimensional examples, we consider a square domain and use $M+1$ grid points in each dimension to discretize the domain. We use $T$ to represent the final time and $N$ is the number of time levels after the initial state (i.e. $\Delta t=T/N$). The computed solution and the exact solution at $t=T$ will be denoted by $u_T$ and $u^*_T$, respectively. The $L_2$ and $L_{\infty}$ norms are computed at $t=T$ defined by $\|u_T-u^*_T\|_2$ and $\|u_T-u^*_T\|_{\infty}$.

\subsection{Heat equations with time-dependent Dirichlet boundary condition}
\label{sec.ex1}
In our first test, we solve the one dimensional heat equation (\ref{eq.heat}). We consider the source $F=\cos(x+t)-\sin(x+t)$ with the initial condition $v=\cos(x)$. The computational domain is $\Omega=[-\pi,\pi]$. The boundary conditions are $f(t)=\cos(t-\pi)$ and $g(t)=\cos(t+\pi)$ so that the exact solution is given as $u=\cos(x+t)$. First we fix our grid size, using $M=100$, to check the accuracy order with respect to $\Delta t$. We set $T=2$. The errors and accuracy order for different $N$ are shown in Table \ref{tab.ex1.t}. Then we fix the time step by using $N=10^5$ and vary $\Delta x$ to test the accuracy in space. The final time $T$ is set to be 0.5. The errors and accuracy orders are shown in Table \ref{tab.ex1.x}. We see from these tables that the accuracy in both time and space are approximately second order defined using both the $L_2$ norm and $L_{\infty}$ norm.

\begin{table}[!ht]
\begin{tabular}{|c|c|c|c|c|}
\hline
 $N$& $\|u_T-u^*_T\|_2$ & rate & $\|u_T-u^*_T\|_{\infty}$ &rate\\
 \hline
 20 & 1.46$\times 10^{-1}$ & & 4.18$\times 10^{-1}$ & \\
 \hline
 40& 4.08$\times 10^{-2}$ & 1.84& 1.18 $\times 10^{-1}$&1.82 \\
 \hline
 80 & 1.05$\times 10^{-2}$ & 1.96 & 3.05$\times 10^{-2}$&1.95 \\
 \hline
 160 & 2.65$\times 10^{-3}$ & 1.99 & 7.72 $\times 10^{-3}$ & 1.98\\
 \hline
\end{tabular}
\caption{(Example \ref{sec.ex1}, one dimensional case) To test the convergence in time, we take $T=2$ and $M=100$. The table shows the $L_2$- and $L_{\infty}$-errors using different $N$.}
\label{tab.ex1.t}
\end{table}

\begin{table}[!ht]
\begin{tabular}{|c|c|c|c|c|}
\hline
 $M$& $\|u_T-u^*_T\|_2$ & rate & $\|u_T-u^*_T\|_{\infty}$ &rate\\
 \hline
 10 & 1.85$\times 10^{-2}$ & & 1.57$\times 10^{-2}$ & \\
 \hline
 20& 2.88$\times 10^{-3}$ & 2.68& 3.53 $\times 10^{-3}$&2.15 \\
 \hline
 40 & 4.78$\times 10^{-4}$ & 2.59 & 8.42$\times 10^{-4}$&2.07 \\
 \hline
 80 & 8.19$\times 10^{-5}$ & 2.55 & 2.05 $\times 10^{-4}$ & 2.04\\
 \hline
\end{tabular}
\caption{(Example \ref{sec.ex1}, one dimensional case) To test the convergence in space, we take $T=0.5$ and $N=10^5$. The table shows the $L_2$- and $L_{\infty}$-errors using different $M$.}
\label{tab.ex1.x}
\end{table}

\begin{table}[!ht]
\begin{tabular}{|c|c|c|c|c|}
\hline
 $N$& $\|u_T-u^*_T\|_2$ & rate & $\|u_T-u^*_T\|_{\infty}$ &rate\\
 \hline
 40 & 1.12$\times 10^{-1}$ &  & 7.13$\times 10^{-2}$& \\
 \hline
 80 & 3.77$\times 10^{-2}$ & 1.57 & 2.38$\times 10^{-2}$&1.58 \\
 \hline
 160 & 1.05$\times 10^{-2}$ & 1.84 & 6.57 $\times 10^{-3}$ & 1.86\\
 \hline
 320 & 2.59$\times 10^{-3}$ & 2.02 & 1.63 $\times 10^{-3}$ & 2.01\\
 \hline
\end{tabular}
\caption{(Example \ref{sec.ex1}, two dimensional case) To test the convergence in time, we take $T=1$ and $M=50$. The table shows the $L_2$- and $L_{\infty}$-errors using different $N$.}
\label{tab.ex2.t}
\end{table}

\begin{table}[!ht]
\begin{tabular}{|c|c|c|c|c|}
\hline
 $M$& $\|u_T-u^*_T\|_2$ & rate & $\|u_T-u^*_T\|_{\infty}$ &rate\\
 \hline
 10 & 5.20$\times 10^{-3}$ & & 3.09$\times 10^{-3}$ & \\
 \hline
 20& 1.18$\times 10^{-3}$ & 2.14& 7.22 $\times 10^{-4}$&2.10 \\
 \hline
 40 & 2.81$\times 10^{-4}$ & 2.07 & 1.73$\times 10^{-4}$&2.06 \\
 \hline
 80 & 6.83$\times 10^{-5}$ & 2.04 & 4.21 $\times 10^{-5}$ & 2.04\\
 \hline
\end{tabular}
\caption{(Example \ref{sec.ex1}, two dimensional case) To test the convergence in space, we take $T=1$ and $N=10^5$. The table shows the $L_2$- and $L_{\infty}$-errors using different $M$.}
\label{tab.ex2.x}
\end{table}

Now we turn to a higher dimensional case where we solve the two dimensional heat equation with the source term given by $b(x,y,t)=2(t+t^2)\sin(x+y)$, the initial condition $r(x,y)=0$ and the time dependent boundary condition $h(x,y,t)=t^2\sin(x+y)$. The exact solution is given by $u=t^2\sin(x+y)$. We set the final time $T=1$. First we fix $M=50$ to test the accuracy order in time. The $L_2$ and $L_{\infty}$ errors are shown in Table \ref{tab.ex2.t}. We can see that the method is almost second order accurate in time. Then we fix $N=10^5$ and test the convergence in space. These results are shown in Table \ref{tab.ex2.x}. We observe also that the method gives approximately second order accurate solution in space.

\begin{table}[!ht]
\begin{tabular}{|c|c|c|c|c|}
\hline
 $N$& $\|u_T-u^*_T\|_2$ & rate & $\|u_T-u^*_T\|_{\infty}$ &rate\\
 \hline
 500 & 3.22$\times 10^{-2}$ & & 1.42$\times 10^{0}$ & \\
 \hline
 1000& 1.73$\times 10^{-2}$ & 0.90& 7.62 $\times 10^{-1}$&0.90 \\
 \hline
 2000 & 5.29$\times 10^{-3}$ & 1.71 & 2.34$\times 10^{-1}$&1.70 \\
 \hline
 4000 & 1.42$\times 10^{-3}$ & 1.90 & 6.36 $\times 10^{-2}$ & 1.88\\
 \hline
\end{tabular}
\caption{(Example \ref{sec.ex3}) To test the convergence in time, we take $T=2$ and $M=1000$. The table shows the $L_2$- and $L_{\infty}$-errors using different $N$.}
\label{tab.ex3.t}
\end{table}
\begin{table}[!ht]
\begin{tabular}{|c|c|c|c|c|}
\hline
 $M$& $\|u_T-u^*_T\|_2$ & rate & $\|u_T-u^*_T\|_{\infty}$ &rate\\
 \hline
 10 & 1.14$\times 10^{-1}$ & & 5.13$\times 10^{-1}$ & \\
 \hline
 20& 3.48$\times 10^{-2}$ & 1.71& 2.46 $\times 10^{-1}$&1.06 \\
 \hline
 40 & 1.15$\times 10^{-2}$ & 1.60 & 1.21$\times 10^{-1}$&1.02 \\
 \hline
 80 & 3.92$\times 10^{-3}$ & 1.55 & 5.97 $\times 10^{-2}$ & 1.02\\
 \hline
\end{tabular}
\caption{(Example \ref{sec.ex3}) To test the convergence in space, we take $T=2$ and $N=10^5$. The table shows the $L_2$- and $L_{\infty}$-errors using different $M$.}
\label{tab.ex3.x}
\end{table}

\subsection{The one dimensional heat equation with the Neumann boundary condition}
\label{sec.ex3}
In this example, we apply the ADE scheme to a problem with the Neumann boundary condition. Consider the following problem
$$
\begin{cases}
u_t=\Delta u +\cos(\pi x)+\pi^2t\cos(\pi x) \, ,\\
u(x,0)=0 \, \mbox{ and } \, u_x(0)=u_x(1)=0 \, .
\end{cases}
$$
The exact solution to this problem is given by $u=t\cos(\pi x)$. We consider the computational domain $\Omega=[0,1]$ and the final time $T=2$. First, we fix $M=1000$ and test the convergence in time using different time steps and $N$. The error in the numerical solutions are shown in Table \ref{tab.ex3.t}.
We also fix $N=10^5$ and vary $M$ to see the convergence in the spatial direction. The results are shown in Table \ref{tab.ex3.x}.

\begin{table}[!ht]
  \centering
  \begin{tabular}{|c|c|c|c|c|}
    \hline
    N & $\|u_T-u^*_T\|_2$& rate & $\|u_T-u^*_T\|_{\infty}$ & rate \\
    \hline
    10 & $5.60\times 10^{-2}$ & & $2.58 \times 10^{-2}$ & \\
    \hline
    20 & $1.66\times 10^{-2}$ &1.76 & $7.65 \times 10^{-3}$ &1.76 \\
    \hline
    40 & $4.50\times 10^{-3}$ & 1.88& $2.08 \times 10^{-3}$ &1.88 \\
    \hline
    80 & $1.17\times 10^{-3}$ &1.94 & $5.40 \times 10^{-4}$ &1.95 \\
    \hline
  \end{tabular}
  \caption{(Example \ref{sec.ex4}) To test the convergence in time, we take $T=0.5$, $M=100$ and $J=200$. The table shows the $L_2$- and $L_{\infty}$-errors using different $N$.}
  \label{tab.ex4.t}
\end{table}

\subsection{The two dimensional time-distributed super-diffusive equation with the time-dependent Dirichlet boundary condition}
\label{sec.ex4}

In this example, we solve the time-distributed super-diffusive equation (\ref{eq.fde}) using the formulation (\ref{eq.ade.scheme1}) by the ADE method. We consider the computation domain $\Omega=(0,\pi)^2$. The functions are chosen as
\begin{eqnarray*}
&& w(\gamma)=\Gamma(7-\gamma) \, ,\, F(x,y,t)=128t^4\sin(x+y)\left[ \frac{360(t-1)}{\ln t}+t^2\right] \, ,\\
&& \psi(x,y,t)=64 t^6 \sin(x+y)
\end{eqnarray*}
so that the exact solution to the problem can be analytically found and is given by $u^*(x,y,t)=\psi(x,y,t)$. First we fix $M=100$ and $J=200$ to test the accuracy order in time. We set $T=0.5$. The result is shown in Table \ref{tab.ex4.t}. The $L_2$ and $L_{\infty}$ norm are both almost second order.

\begin{figure}[!ht]
(a)\includegraphics[width=0.45\textwidth]{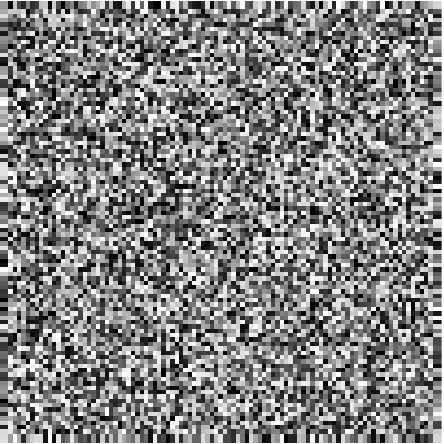}
(b)\includegraphics[width=0.45\textwidth]{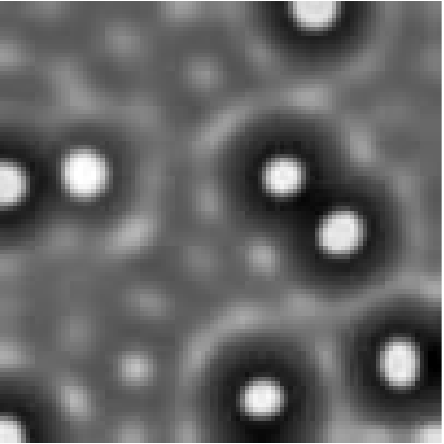}
(c)\includegraphics[width=0.45\textwidth]{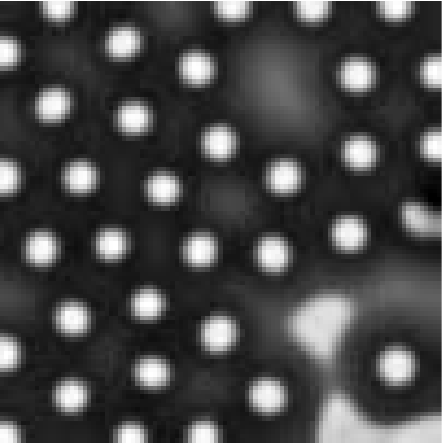}
(d)\includegraphics[width=0.45\textwidth]{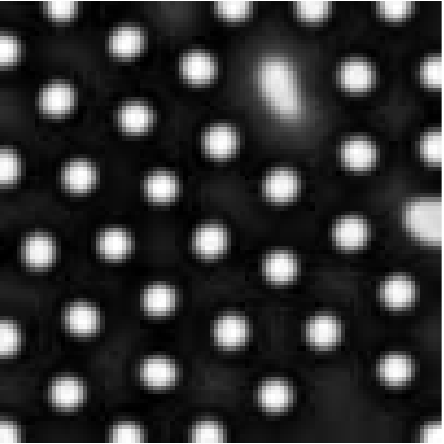}
(e)\includegraphics[width=0.45\textwidth]{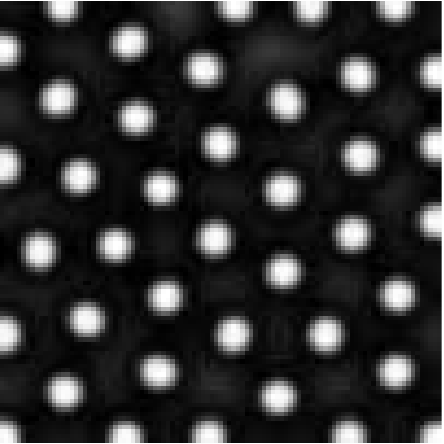}
(f)\includegraphics[width=0.45\textwidth]{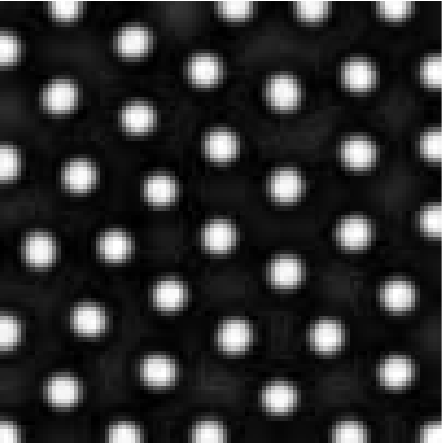}
\caption{(Example \ref{sec.ex6}) The initial condition is shown in (a). Turing pattern obtained after (b) 2000 iterations, (c) 6000 iterations, (d) 12000 iterations, (e) 16000 iterations and (f) 20000 iterations.}
\label{fig.turing}
\end{figure}

\subsection{The two dimensional sub-diffusive activator-inhibitor system with the zero Neumann boundary condition}
\label{sec.ex6}

In this example, we use the ADE scheme to solve the sub-diffusive activator-inhibitor system using $100\times100$ grid points with $\Delta x=1$. The coefficient is chosen as $a_{11}=0.899$, $a_{22}=-0.91$, $a_{12}=1$, $a_{21}=-a_{11}$, $\delta=2$, $D=0.516$, $\alpha=0.92$ and $\beta=0.88$. With this setting, there is only one fixed point for this system given by $u=v=0$. For the initial condition, we add some uniformly distributed noise ranging from -0.1 to 0.1. In equation (\ref{eq.approxFD}) when approximating the fractional derivative, we find that the number of terms in the expression grows rapidly in time yet the magnitude of some of these terms are rather small. To improve the computational efficiency, we follow the approach in \cite{hernandez2009dynamics} and keep only the latest 800 terms whose binomial coefficient $\binom{\alpha}{m}$ or $\binom{\beta}{m}$ is larger than $10^{-7}$. In our implementation, we use the time step $\Delta t=0.08$. Our result is shown in Figure \ref{fig.turing}. The method as stated in \cite{hernandez2009dynamics} requires approximately $3.8\times 10^5$ iterations to reach the steady state, while our ADE scheme takes only $2\times 10^4$ iterations (roughly $1/20$ of the former method) to obtain the final pattern since the unconditionally stable property in the ADE method allows us to use a significantly larger time step in updating the solution.

\section{Conclusion}
\label{sec.conclusion}

The ADE scheme is an unconditionally stable fully explicit AOS scheme and can efficiently solve time-evolution equations. We have extended the ADE scheme to time evolution equations with the time dependent Dirichlet boundary condition and the Neumann boundary condition. We also present various numerical examples including the heat equation with the time dependent Dirichlet boundary condition and the Neumann boundary condition, and two applications involving fractional derivatives including the time-distributed order super-diffusive equation and a reaction-diffusion system in the sub-diffusive regime. These numerical tests have demonstrated that the numerical approach is simple to implement and is computational efficient for a wide range of applications.

\section*{Acknowledgment}
The work of Leung was supported in part by the Hong Kong RGC grants 16303114 and 16309316.

\bibliographystyle{plain}
\bibliography{liu,syleung}

\end{document}